\documentclass[a4paper,11pt]{article}
\usepackage{mathrsfs}
\usepackage{latexsym,epsfig,amsmath,bm}
\usepackage{amssymb}
\usepackage{rawfonts}

\begin{document}

\input prepictex
\input pictex
\input postpictex

\newcommand\bp{\beginpicture}
\newcommand\ep{\endpicture}

\newcounter{prostredi}
\def\theprostredi{\arabic{prostredi}}

\parindent=20pt
\makeatletter
\def\@begintheorem#1#2{\trivlist
   \item[\hskip \labelsep{\bfseries #1\ #2.}]} %\itshape}
\def\@opargbegintheorem#1#2#3{\trivlist
      \item[\hskip \labelsep{\bfseries #1\ #2.\ (#3.)}]} % \itshape}
\def\@endtheorem{\endtrivlist}
\makeatother

\makeatletter \@addtoreset{equation}{section}
\def\theequation{\thesection.\arabic{equation}}
\makeatother

\newenvironment{claim}{\par\bigskip\noindent%
\refstepcounter{prostredi}{\bf Claim \theprostredi.}\quad\bgroup\sl
}
{\egroup\par\bigskip\endtrivlist}%

\newcommand{\forw}{\hspace{-\mathsurround}\rlap{\raisebox{1.6ex}{%
        \vphantom{-}\smash{$\rightarrow$}}}\hspace{\mathsurround}}
\newcommand{\back}{\hspace{-\mathsurround}\rlap{\raisebox{1.6ex}{%
        \vphantom{-}\smash{$\leftarrow$}}}\hspace{\mathsurround}}
\newcommand{\deltaback}{\back \Delta}
\newcommand{\deltaforw}{\forw \Delta}
\newcommand{\cforw}{\forw C}
\newcommand{\vuforw}{\forw vu}
\newcommand{\vwforw}{\forw vw}

\newcommand{\qed}{$\qquad\square$\bigbreak}

\def\R{{\mathbb R}}
\def\Z{{\mathbb Z}}
\def\N{{\mathbb N}}
\newcommand{\diam}{{\mbox{diam}}}

\newcounter{vetyc}
\def\thevetyc{\Alph{vetyc}}
\newcommand{\mezera}{\vspace{3.34mm}}
\newenvironment{vetacite}[1]{\par\mezera\noindent%
\refstepcounter{vetyc}{\bf Theorem
\thevetyc{} {\cite{#1}}.}\quad\bgroup\sl }
{\egroup\par\mezera\endtrivlist}%

\newcommand\bvec{\begin{vetacite}}
\newcommand\evec{\end{vetacite}}

\newtheorem{theorem}{Theorem}%[section]
\newtheorem{lemma}[theorem]{Lemma}
\newtheorem{definition}[theorem]{Definition}
\newtheorem{conjecture}[theorem]{Conjecture}
\newtheorem{corollary}[theorem]{Corollary}
\newtheorem{proposition}[theorem]{Proposition}
\newtheorem{observation}[theorem]{Observation}
\newcommand{\cl}{{\rm cl}}

\renewcommand{\theenumi}{\arabic{enumi}}
\renewcommand{\labelenumi}{(\theenumi)}

\newcommand{\bd}[2] {\partial_{#1}(#2)}
\newcommand{\Setx}[1] {\left\{#1\right\}}
\newcommand{\size}[1] {\left|{#1}\right|}
\newcommand{\symm} {\triangle}
\begin{center}
\LARGE Star subdivisions and connected even factors in the square of a graph  

\vspace{6mm}

\large Jan Ekstein $^1$,  P\v{r}emysl Holub $^1$, Tom\'a\v s Kaiser \footnote {Department of Mathematics, University of West Bohemia, and Institute for Theoretical Computer Science
(ITI), Charles University, Univerzitni 22, 306 14 Pilsen, Czech Republic, e-mail: $\{$ekstein, holubpre, kaisert$\}$@kma.zcu.cz; research supported by Grant No. 1M0545 of the Czech Ministry of Education.}$^,$\footnote{Research supported by Research Plan MSM 4977751301 of the Czech Ministry of Education and the grant GA\v CR 201/09/0197 of the Czech Science Foundation.}, \\
 Liming Xiong \footnote{  Department of Mathematics, Beijing Institute of Technology, Beijing, 100081, P.R. China; Supported by Nature Science Funds of China.},
 Shenggui Zhang \footnote{ Department of Applied Mathematics, Northwestern Polytechnical University, Xian, Shaanxi 710072, P.R. China; Supported by Nature Science
Funds of China.}

\end{center}

\vspace{5mm}

\centerline{\bf Abstract.}

\vspace{2mm}

\small
\hspace*{1.2cm}\begin{minipage}[l]{0.85\textwidth}
For any positive integer $s$, a $[2,2s]$-factor in a graph $G$ is a connected even factor with maximum degree at most $2s$. We prove that if every induced $S(K_{1, 2s+1})$ in a graph $G$ has at least $3$ edges in a block of degree at most two, then $G^2$ has a $[2,2s]$-factor. This extends the results of Hendry and Vogler and of Abderrezzak et al.
\vspace{5mm}


\vspace{5mm}


\end{minipage}

\normalsize

\vspace{3mm}

\section{Introduction}

We use Bondy and Murty \cite{bon} for
terminology and notation not defined here and we consider only finite
undirected simple graphs, unless otherwise stated.

Let $G=(V,E)$ be a graph with  vertex set $V$ and  edge
set $E$. Let $\alpha(G)$ denote the independence number of $G$, i.e., the cardinality of a largest independence set in $G$. For any vertex $x$ of $G$, let $d_G(x)$ denote the degree of $x$ in $G$, $N_G(x)$ the set of all neighbors of $x$ in $G$, $N_G[x]=N_G(x)\cup \{x\}$. 
The square of a graph $G$, denoted by $G^2$, is the graph with $V(G^2)=V(G)$ in which two vertices are adjacent if their distance in $G$ is at most two. Thus $G\subseteq G^2$. 

 For any $S\subseteq V(G)$, we denote by $G[S]$ the subgraph of $G$
induced by $S$. 
For a positive integer $s$, the graph  $S(K_{1,2s+1})$ is obtained from the complete bipartite graph $K_{1,2s+1}$ by subdividing each edge once. The graph $G$ is said to be {\sl $S(K_{1,2s+1})$-free} if it does not contain any induced copy of $S(K_{1,2s+1})$.

A connected graph that has no cut vertices is called a {\em block}.
{\em A block of a graph} $G$ is a subgraph of $G$ that is a block and is maximal with respect to this property. {\sl The degree of a block} $B$ in a graph $G$, denoted by $d(B)$, is the number of cut vertices of $G$ belonging to $V(B)$.

A {\sl factor } in a graph $G$ is a spanning subgraph of $G$. A connected even factor in $G$ is a connected factor in $G$ with all vertices of even degree.
A {\sl $[2,2s]$-factor} in $G$ is a connected even factor in $G$ in which degree of every vertex is at most $2s$. A graph is {\sl hamiltonian} if it has a  spanning cycle.
 In other word, a graph is hamiltonian if and only if it has a $[2,2]$-factor.

 The following result concerns the existence of a $[2,2]$-factor in the square of a 2-connected graph.

\bvec{F1} \label{hamsquare}
Let $G$ be a 2-connected graph. Then $G^2$ is hamiltonian.
\evec

Gould and Jacobson in \cite{Gould} conjectured that for the hamiltonicity of $G^2$, the connectivity condition can be relaxed for $S(K_{1,3})$-free graphs. Their conjecture was proved by Hendry and Vogler in \cite{Hendry}.

\bvec{Hendry} \label{Hendry} Let $G$ be a connected $S(K_{1,3})$-free graph. Then $G^2$ is hamiltonian, $i.e.,$ has a $[2,2]$-factor.
\evec

Moreover, Abderrezzak, Flandrin and Ryj\'a\v cek in \cite{Ryjacek} proved the following result in which graphs may contain an induced $S(K_{1,3})$ of a special type.

\bvec{Ryjacek} \label{Ryjacek} Let $G$ be a connected graph such that every induced $S(K_{1,3})$ in $G$ has at least three edges in a block of degree at most two. Then $G^2$ is hamiltonian, $i.e.,$  has a $[2,2]$-factor.
\evec

It is a natural question if there exists a $[2,2s]$-factor in the square of a graph if one replaces $S(K_{1,3})$ by $S(K_{1,2s+1})$ in Theorems~\ref{Hendry} and \ref{Ryjacek}. In this paper, we will give a positive  answer to this question; we will extend Theorems~\ref{Hendry} and \ref{Ryjacek} as follows.

\begin{theorem}
\label{square} Let $G$ be a  connected $S(K_{1,2s+1})$-free graph of order at least three and $s$ a positive integer. Then $G^2$ has a   $[2,2s]$-factor.
\end{theorem}

Since the square of an $S(K_{1,2s+1})$ itself has no $[2,2s]$-factor, Theorem~\ref{square} is the best possible in a sense.

\begin{theorem}
\label{thinduced} Let $s$ be a positive integer and $G$ be a connected graph such that every induced $S(K_{1, 2s+1})$ has at least three edges in a block of degree at most two. Then $G^2$ has a $[2, 2s]$-factor.
\end{theorem}

Note that Theorem \ref{thinduced} is a strengthening of Theorem \ref{square}, but we state Theorem \ref{square} separately because it will be used in the proof of Theorem \ref{thinduced}.

\section{Preliminaries and auxiliary results}\label{sec:auxi}
As noted in Section 1, for graph-theoretic
notation not explained in this paper, we refer the reader to
\cite{bon}.

A graph $G$ is {\sl even} if every vertex of $G$ has even degree.
In the subsequent sections, we frequently take the symmetric
difference of two subgraphs of a graph. Let $H,\, H'$ be subgraphs of a graph
$G$. The graph $H\vartriangle H'$ has vertex set $V(H)\cup V(H')$ and its edge set is the symmetric difference of $E(H)$ and $E(H')$.
 Note that if $H$ and $H'$ are both even
graphs, then $H\vartriangle H'$ is also an even graph.

A {\sl trail} between vertices $u_0$ and $u_r$ is a finite sequence $T=u_0e_1u_1e_2u_2\cdots
e_ru_r$, whose terms are alternately vertices and edges, with
$e_i=u_{i-1}u_i, \, 1\le i\le r$, where the edges are distinct. A
trail $T$ is \textit{closed} if $u_0=u_r$, and it is
\textit{spanning} if $V(T)=V(G)$. An $s$-trail between $u_0$ and $u_r$ is a trail starting at $u_0$, ending at $u_r$ and in which every vertex is visited at most $s$ times. In other words, a $[2,2s]$-factor in a graph $G$ can be viewed as a spanning closed $s$-trail in $G$ and vice versa. We define the degree of a vertex $x$ in an $s$-trail as the number of edges incident with $x$ in the corresponding $[2,2s]$-factor.

We use the following fact (see \cite{jackson}, Corollary 2.3.1 for a proof).

\bvec{jackson} \label{ordaz}
Let $k\geq 2$ be an integer and $G$ a $k$-connected graph. If $\alpha(G)>k$ then $V(G)$ can be covered with $\alpha(G)-k$ disjoint paths.
\evec

From the proof of this Theorem it follows that the statement is true without the restrictions on $k$, in particular for $k=0$.

 \begin{corollary}\label{cor1}
 Let $G$ be a graph. Then there are at most $\alpha (G)$ disjoint paths covering $V(G)$.
 \end{corollary}
 
Let $G_1, G_2$  be graphs such that $V(G_1)\cap V(G_2)=\{x\}$. The symbol $G=G_1xG_2$ denotes a 
graph $G$ with $V(G)=V(G_1)\cup V(G_2)$ and $E(G)=E(G_1)\cup E(G_2)$.

Given a subgraph $K$ of a graph $H$, we define $\bd H K$ as the set of
all edges of $H$ with exactly one endvertex in $V(K)$. Thus $\bd H K$
is a (not necessarily minimal) edge-cut.

\begin{lemma}\label{lemtrail} Let $H$ be a connected graph and $P= xyz$ a
  path of length two such that $V(H)\cap V(P)=\{x\}$. If $(HxP)^2$ has a $[2,2s]$-factor, then one of
  the following holds:
 \begin{itemize}
 \item [(a)] $H^2$ contains a spanning closed $s$-trail $T$ such that the degree
   of $x$ in $T$ is at most $2s-2$, or
 \item [(b)] $H^2$ contains a spanning $s$-trail $T$ between $x$ and some $x'\in N_H(x)$.
 \end{itemize}     
\end{lemma}

{\bf Proof.} Let $F$ be a $[2,2s]$-factor of $(HxP)^2$ and let
$K_0,\dots,K_\ell$ be all the components of $F\setminus\Setx{y,z}$,
where $x\in V(K_0)$. Furthermore, define $W = N_F(y)\setminus\Setx z$
and $W_i = W\cap V(K_i)$ ($i=0,\dots,\ell$). Observe that each $W_i$
is nonempty. Clearly, the induced subgraph $Q$ of $H^2$ on $W\cup\Setx
x$ is complete.

Since $F$ covers $z$, it includes the edges $yz$ and $xz$. For $0\leq
i \leq \ell$, every edge in $\bd F {K_i}$ is incident with $y$, except
for the edge $xz\in\bd F {K_0}$. Since
\begin{equation*}
  \bd F {K_i} = \bd H {K_i} \cap E(F)
\end{equation*}
and the intersection of any edge-cut with an eulerian subgraph has
even cardinality, we conclude that for $0\leq i \leq \ell$,
\begin{equation*}
  \size{W_i}\text{ is odd if and only if $i=0$.}
\end{equation*}
If $w\in W_i$ and $w\neq x$, then the degree of $w$ in $K_i$ is odd
and does not exceed $2s-1$. The same is true for $w=x$ provided that
$x\notin W$, since then $xz$ is the only edge of $\bd F {K_0}$
incident with $x$. On the other hand, if $x\in W$, then both $xz$ and
$xy$ have this property, so the degree of $x$ in $K_0$ is even and
does not exceed $2s-2$.

For each $i$, $0\leq i \leq\ell$, choose a matching $M_i$ that covers
all except one or two vertices of $W_i$ (one if $i=0$, two otherwise)
and uses as few edges as possible from $F$. We argue that the
symmetric difference $K_i \vartriangle M_i$ is connected. We may assume that
$M_i$ uses at least one edge of $F$, otherwise there is nothing to
prove. For a fixed $i$, let $X \subseteq W_i$ be the set consisting of
vertices incident with edges in $E(M_i)\cap E(F)$, together with the
vertices of $W_i$ left uncovered by $M_i$. By the choice of $M_i$,
$K_i[X]$ must be complete and $\size X \geq 3$. All the edges of $K_i$
that are removed as a result of taking the symmetric difference are
edges of $K_i[X]$. Since any graph obtained by removing a matching
from a complete graph on at least 3 vertices is connected, the claim
follows. 

Observe that for $i\geq 1$, each $K_i \vartriangle M_i$ contains exactly two
vertices of odd degree (and the degree does not exceed $2s-1$). The
same is true for $i=0$ unless $x\in W$ and $x$ is not incident with
$M_0$, in which case $K_0 \vartriangle M_0$ is eulerian and the degree of $x$
in this graph is at most $2s-2$. It follows that if $\ell=0$, then we
can set $T := K_0 \vartriangle M_0$ and we are done ($T$ satisfies condition
(a) if $x\in W\setminus V(M_0)$ and condition (b) otherwise).

If $\ell \geq 1$, then let $u_0$ be the vertex of $W_0\setminus
V(M_0)$, and for $i\geq 1$, let $W_i\setminus V(M_i) =
\Setx{u_i,v_i}$. Taking the union of all the graphs $K_i\vartriangle M_i$ and
adding the edges $u_0v_1, u_1v_2, \dots, u_{\ell-1}v_\ell$, we obtain
a connected graph $T$ in which the only vertices of odd degree are $x$
and $u_\ell$, and which satisfies condition (b) in the lemma. $\blacksquare$

  Using a similar argument as in the proof of Lemma~\ref{lemtrail}, one can prove the following.

\begin{lemma}\label{lemtrail1}Let $H$ be a connected graph and $P= xy$ an edge such that $V(H)\cap V(P)=\{x\}$. If $(HxP)^2$ has a   $[2,2s]$-factor, then
  $H^2$ has a spanning $s$-trail $T$ between $x'\in N_H[x]$ and some vertex $x''\in N_H(x)$.
 \end{lemma}

The following theorem will be used in the proof of Theorem \ref{thinduced}.

 \bvec{F1} \label{thfleischner} Let $y$ and $z$ be arbitrarily chosen vertices of a 2-connected graph $G$. Then $G^2$ has a hamiltonian cycle $C$ such that the edges of $C$ incident with $y$ are in $G$ and at least one of the edges of $C$ incident with $z$ is in $G$. If $y$ and $z$ are adjacent in $G$, then these are three different edges.
 \evec

\section{Proofs}\label{sec2}

The purpose of this section is to prove Theorem \ref{thinduced}. As mentioned in Section 1, the proof makes use of Theorem \ref{square} which we derive next.

{\bf Proof of Theorem~\ref{square}.} This proof is inspired by the proof in \cite{Hendry}.
We prove our result by induction on $|V(G)|$. Clearly $G^2$ is hamiltonian (hence has a $[2,2]$-factor) for graphs with $|V(G)|\leq 6$, since $G$ is $S(K_{1,3})$-free. By Theorem~\ref{hamsquare}, we may assume that $G$ has cut vertices. If all cut vertices have degree two, then $G$ is a path and hence $G^2$ is hamiltonian. So we may assume that there is a cut vertex $u$ such that $d_G(u)=d\geq 3$. Since $G$ is connected,
 we may take a spanning tree $S$ of $G$ such that $S$ contains all edges of $G$ incident with $u$. We label the neighbors of $u$ by $u_1, u_2, \cdots, u_d$ in such a way that $d_G(u_i)\geq 2$ for $1\leq i\leq m$ and $d_G(u_i)=1$ for $m+1\leq i\leq d$. For $i\leq m$, let $G_i$ be the subgraph of $G$ induced by the vertices in the component of the forest $S-u$ containing $u_i$; we fix a neighbour $u_i'$ of $u$ that is not contained in the same component of $G-u$ as $u_i$ (note that there must be such a vertex since $u$ is a cut vertex of $G$), and let $H_i=G[V(G_i)\cup \{u, u_i'\}]$. Then $H_i$ is a proper $S(K_{1, 2s+1})$-free subgraph of $G$ since $H_i$ is an induced subgraph of $G$ and $d_G(u)\geq 3$. Note that $H_i$ is connected. By the inductive hypothesis, $H_i^2$ has a   $[2, 2s]$-factor. Note that $d_{H_i^2}(u_i')=2$.

By Lemma \ref{lemtrail} it follows that at least one of the following facts holds.
\begin{itemize}
\item [(a)] there exists a spanning closed $s$-trail $T_i$ in $G_i^2$ such that $d_{T_i}(u_i)\leq 2s-2$;
\item [(b)] there exists a spanning $s$-trail $T_i$ in $G_i^2$ between $u_i$ and some $z_i\in N_{G_i}(u_i)$ .
\end{itemize}

Without loss of generality we may assume that $\{ u_1, u_2, \dots, \, u_{m'} \} \subseteq \{ u_1, u_2, \dots, \, u_m \}$ is the set of all vertices $u_i$ such that $G_i$ has an $s$-trail of type (b), for a suitable $m'\leq m$.
Construct the graph $H$ from $G[\{ u_1, u_2, \dots , u_{m'}, z_1, z_2, \dots , z_{m'}\}]$ by contracting edges $u_iz_i$ to a vertex $w_i$ for $i=1, \dots, m'$. Since $G$ is $S(K_{1, 2s+1})$-free,  $\alpha (H)\leq 2s$. By Corollary~\ref{cor1}, there are  $\ell\leq \alpha (H)$ vertex-disjoint paths  $P_1, P_2, \dots, P_\ell$ covering $V(H)$. Without loss of generality, we may assume that  $P_i=w_{s_{i-1}+1}w_{s_{i-1}+2}\dots w_{s_i}$, for $ i=1, \dots, \,  \ell$ (where $s_0=0$ and $s_\ell=m'$). Since we contracted edges $u_jz_j$ to vertices $w_j$, both $u_j$ and $z_j$ have a neighbor in $\{ u_{j+1}, z_{j+1}\}$ in $G^2$ for $i=1, \dots, \ell$, and $j=s_{i-1}+1, \dots,  s_i-1$. Hence from the paths $P_i$ ($i=1, \dots, \ell$) and $s$-trails $T_j$ ($i=1, \dots, m'$) we can obtain the following $s$-trails $F_i$ in $G^2$:
\begin{itemize}
\item[-]  for a trivial (one-vertex) path $P_i$, $F_i=T_i$,
\item[-]  for a nontrivial path $P_i$, we construct $F_i$ by joining the trails $T_{s_{i-1}+1}, \dots, \, T_{s_{i}}$ with the edges $x_jx_{j+1}$, where $x_j\in \{u_j,z_j\}$ and $x_{j+1}\in \{ u_{j+1},z_{j+1} \}$ with respect to $P_i$. Clearly $d_{F_i}(u_{s_{i-1}+1})<2s$, $d_{F_i}(x_{s_i})<2s$ 
and $F_i$ spans all the vertices of $G_{s_{i-1}+1}\cup \dots \cup G_{s_{i}}$.
\end{itemize}
Note that the number of $s$-trails $F_i$ is $\ell\leq 2s$.

Let $T=u_{m'+1}T_{m'+1}u_{m'+1}u_{m'+2}T_{m'+2}u_{m'+2}\dots u_m T_m u_m u_{m+1} \dots u_d$ be an $s$-trail containing all vertices of $G_{m'+1} \cup \dots \cup G_m$ and all neighbours of $u$ of degree one in $G$. We set $F'= u_1 F_1x_{s_1}ux_{s_2}F_2u_{s_1+1}u_{s_2+1}F_3 \dots x_{s_\ell}F_\ell u_{s_{\ell-1}+1}u_{m'+1}$ for even $\ell$ and $F'=u_1  F_1x_{s_1}ux_{s_2}F_2u_{s_1+1}u_{s_2+1}F_3 \dots u_{s_{\ell-1}+1}F_\ell x_{s_\ell}u u_{m'+1}$ for odd $\ell$. In both cases, $F'$ is an $s$-trail containing all vertices of $G_1 \cup \dots \cup G_{m'}$. Finally, we construct a trail $F=u_1F'u_{m'+1}Tu_du_1$.
Clearly, $d_F(u)=\ell\leq 2s$ and $F$ corresponds to a $[2,2s]$-factor in $G^2$.
\hfill $\blacksquare$

\begin{corollary}Let $G$ be a simple connected graph with $\Delta (G)\leq 2s$. Then $G^2$ has a   $[2,2s]$-factor.
\end{corollary}

 Before we present the proof  of Theorem~\ref{thinduced},  we give some additional definitions. Let $x$ be a cut vertex of $G$, and $H'$ be a component of $G-x$. Then the subgraph $H=G[V(H')\cup \{x\}]$ is called a {\em branch} of $G$ at $x$. Let $F$ be a connected subgraph of $G$ and $x$ some vertex of $F$. Let $P_i(x)$  denote a path on $i$ vertices with end vertex $x$. The subgraph $F$ is called to be {\em nontrivial} at $x$ if it contains a $P_3(x)$ as a proper induced subgraph ($i.e.,$ $F$ is trivial at $x$ if $F=P_3(x)$ or $V(F)\subseteq N[x]$).

Now we present the proof of Theorem~\ref{thinduced}.

{\bf Proof of Theorem~\ref{thinduced}.} We prove this theorem by contradiction. Suppose that Theorem~\ref{thinduced} is not true and choose a graph $G$ in such a way that
\begin{itemize}
\item [(1)] $G$ is connected and every induced $S(K_{1, 2s+1})$ in $G$  has at least three edges in a block of degree at most two;
\item [(2)] $G^2$ has no $[2, 2s]$-factor;
\item [(3)] $|V(G)|$ is minimized with respect to (1) and (2).
\end{itemize}

The following fact is necessary for our proof.

\begin{claim}\label{claim11} Let $x$ be a cut vertex of $G$ and $F_1,F_2$ two connected subgraphs of $G$ such that $F_1$, $F_2$ belong to different branches of $G$ at $x$. Assume that $F_2$ is nontrivial at $x$, $i.e.,$ $F_2$ contains an induced $P_3(x)=xyz$ as a proper induced subgraph. Then the graph $G'=F_1xP_3(x)$ also satisfies (1).
\end{claim}

{\it Proof of Claim~\ref{claim11}.} If not, there exists in $G'$ some $S(K_{1,2s+1})$ that has no connected part of order at least $4$ in a block of degree at most two. But if so, it is the same in $G$, since any $S(K_{1,2s+1})$ in $G'$ is also an induced $S(K_{1,2s+1})$ of $G$. \hfill $\Box$

Since in our proof we have assumed that $G^{2}$ has no $[2,2s]$-factor, we know from Theorem 1 that $G$
 contains some $S(K_{1,2s+1})$ as an induced subgraph. By (1), the $S(K_{1,2s+1})$ has at least 3 edges in some block $H$
 of $G$ of degree at most 2. Notice that $|V(H)| \geq 5$.

 \emph{Case 1:} $d(H) = 1$. Let $c$ be the cut vertex of $G$ belonging to $H$ and let $R$ be the union of all branches
 of $G$ at $c$ which intersect $H$ only at $c$.

 If $H$ is trivial at $c$, then $V(H) - \{c\} = \{b_{1}, b_{2}, ..., b_{h}\} \subseteq N(c)$. The graph
 $G' = Rc(cb_{1})$ satisfies condition (1). So by minimality of $G$, the graph $G'^{2}$ has a $[2, 2s]$-factor and,
 by Lemma 5, $R^{2}$ has a spanning $s$-trail $T$ between some $c' \in N_{R}[c]$ and some $c'' \in N_{R}(c)$. Let
 $F = c'Tc''b_{1}...b_{h}c'$. It is easy to see that $F$ is a $[2,2s]$-factor in $G^{2}$, a contradiction.

 Hence $H$ is nontrivial at $c$, i.e., it contains a proper induced path $P_{3}(c) = cb_{1}b_{2}$. By Theorem E,
 $H^{2}$ contains a hamiltonian path $b_{1}P_{H^{2}}c$ connecting $b_{1}$ and $c$. On the other hand the graph
 $G'' = RcP_{3}(c)$ is connected and, by Claim \ref{claim11}, $G''$ satisfies condition (1). Since
 $|V(G'')| < |V(G)|$, $(G'')^{2}$ has a $[2, 2s]$-factor and by Lemma 4, one of the following subcases occur.

 If the graph $R^{2}$ has a spanning closed $s$-trail $T'$ in which $d_{T'}(c) \leq 2s - 2$, then $F = cT'cb_{1}P_{H^{2}}c$
 is a $[2, 2s]$-factor in $G^{2}$, a contradiction.

 If the graph $R^{2}$ has a spanning $s$-trail $T''$ between $c$ and some neighbor $c''' \in N_{R}(c)$, then
 $F = cT''c'''b_{1}P_{H^{2}}c$ is a $[2, 2s]$-factor in $G^{2}$, contradicting condition (2).

 \emph{Case 2:} $d(H) = 2$. Let $c_{1}$ and $c_{2}$ be two cut vertices of $G$ belonging to $H$ and let $B_{i}$,
 $i = 1, 2$, be the union of all branches of $G$ at $c_{i}$ not containing $H$. This means that
 $G = (B_{1}c_{1}H)c_{2}B_{2}$. The subgraph $H$ is a block and thus, by Theorem E, $V(H)$ can be covered by two vertex-disjoint paths $a_{1}P_{H}^{1}a_{2}$ and $c_{2}P_{H}^{2}c_{1}$ in $H^{2}$, where $a_{1} \in N(c_{1})$ and
 $a_{2} \in N(c_{2})$. We distinguish, up to symmetry, the following three subcases.

 \emph{Subcase 2.1:} $B_{1}$ is trivial at $c_{1}$ and $B_{2}$ is trivial at $c_{2}$.

 If $V(B_{1}) = \{b_{1}, b_{2}, ..., b_{k}, c_{1}\} \subseteq N[c_{1}], k \geq 1$, and $B_{2} = P_{3}(c_{2}) =
 c_{2}d_{1}d_{2}$, then $F = c_{1}b_{1}b_{2}...b_{k}a_{1}P_{H}^{1}a_{2}d_{1}d_{2}c_{2}P_{H}^{2}c_{1}$ is even
 a hamiltonian cycle in $G^{2}$, which contradicts the fact that $G^2$ has no $[2,2s]$-factor.

 The proof is similar if $B_{1} = P_{3}(c_{1})$ and $V(B_{2}) \subseteq N[c_{2}]$.

 If $V(B_{1}) = \{b_{1}, b_{2}, ..., b_{k}, c_{1}\} \subseteq N[c_{1}]$ and $V(B_{2}) = \{d_{1}, d_{2}, ..., d_{l},
 c_{2}\} \subseteq N[c_{2}]$, then $F = c_{1}b_{1}b_{2}...b_{k}a_{1}P_{H}^{1}a_{2}d_{1}d_{2}...d_{l}c_{2}P_{H}^{2}c_{1}$
 is also a hamiltonian cycle in $G^{2}$, contradicting (2).

 Finally, if $B_{1} = P_{3}(c_{1}) = c_{1}b_{1}b_{2}$ and $B_{2} = P_{3}(c_{2}) = c_{2}d_{1}d_{2}$, then again the cycle
 $F = c_{1}b_{2}b_{1}a_{1}P_{H}^{1}a_{2}d_{1}d_{2}c_{2}P_{H}^{2}c_{1}$ gives a similar contradiction.

 \emph{Subcase 2.2:} $B_{1}$ is nontrivial at $c_{1}$ and $B_{2}$ is trivial at $c_{2}$.

 Since $|V(H) \cup V(B_{2})| > 3$, there exists some vertex in $V(H) \cup V(B_{2})$ (for example each vertex in
 $V(B_{2}) \setminus \{c_{2}\}$) nonadjacent to $c_{1}$, the subgraph $G' = Hc_{2}B_{2}$ is nontrivial. Then $G'$
 contains a path $P_{3}(c_{1}) = c_{1}n_{1}n_{2}$ as a proper induced subgraph. Now let $G_{1} = B_{1}c_{1}n_{1}n_{2}$.
 By Claim \ref{claim11}, $G_{1}$ satisfies condition (1). By minimality of $G$, the graph $G_{1}^{2}$ has
 a $[2, 2s]$-factor and thus, by Lemma 4, we have the following two possibilities.

 $a)$ The graph $B_{1}^{2}$ has a spanning closed $s$-trail $T$ in which $d_{T}(c_{1}) \leq 2s - 2$.

     If $V(B_{2}) = \{b_{1}, b_{2}, ..., b_{k}, c_{2}\} \subseteq N[c_{2}], k \geq 1$, then $F =
     c_{1}Tc_{1}a_{1}P_{H}^{1}a_{2}b_{1}b_{2}...b_{k}c_{2}P_{H}^{2}c_{1}$ is a $[2, 2s]$-factor in $G^{2}$, a contradiction with (2).

     If $B_{2} = P_{3}(c_{2}) = c_{2}d_{1}d_{2}$, then $F = c_{1}Tc_{1}a_{1}P_{H}^{1}a_{2}d_{1}d_{2}c_{2}P_{H}^{2}c_{1}$
     is a $[2, 2s]$-factor in $G^{2}$, which contradicts condition (2).

 $b)$ The graph $B_{1}^{2}$ has a spanning $s$-trail $T'$ between $c_{1}$ and some neighbor $c_{1}' \in
      N_{B_{1}}(c_{1})$.

     If $V(B_{2}) = \{b_{1}, b_{2}, ..., b_{k}, c_{2}\} \subseteq N[c_{2}], k \geq 1$, then $F =
     c_{1}T'c_{1}'a_{1}P_{H}^{1}a_{2}b_{1}b_{2}...b_{k}c_{2}P_{H}^{2}c_{1}$ is a $[2, 2s]$-factor in $G^{2}$
     and contradicts (2).

     If $B_{2} = P_{3}(c_{2}) = c_{2}d_{1}d_{2}$, then $F = c_{1}T'c_{1}'a_{1}P_{H}^{1}a_{2}d_{1}d_{2}c_{2}
     P_{H}^{2}c_{1}$ is a $[2, 2s]$-factor in $G^{2}$, a contradiction with (2).

 \emph{Subcase 2.3:} $B_{1}$ is nontrivial at $c_{1}$ and $B_{2}$ is nontrivial at $c_{2}$.

 Let $G_{1}$ be the same graph as in Subcase 2.2 and in a similar way as in Subcase 2.2 let $G_{2} =
 B_{2}c_{2}m_{1}m_{2}$, where a path $c_{2}m_{1}m_{2}$ is a proper induced subgraph of $Hc_{1}B_{1}$. Then, by Claim \ref{claim11}, 
 both $G_{1}$ and $G_{2}$ satisfy condition (1). By minimality of $G$, the graphs $G_{1}^{2}$ and
 $G_{2}^{2}$
 have a $[2, 2s]$-factor and thus, by Lemma 4, we have the following two possibilities.

 $a)$ The graph $B_{1}^{2}$ has a spanning closed $s$-trail $T$ in which $d_{T}(c_{1}) \leq 2s - 2$.

      If the graph $B_{2}^{2}$ has a spanning closed $s$-trail $T'$ in which $d_{T'}(c_{2}) \leq 2s - 2$, then
      $F = c_{1}Tc_{1}a_{1}P_{H}^{1}a_{2}c_{2}T'c_{2}P_{H}^{2}c_{1}$ is a $[2, 2s]$-factor in $G^{2}$
      and contradicts (2).

      If the graph $B_{2}^{2}$ has a spanning $s$-trail $T''$ between $c_{2}$ and some neighbor $c_{2}' \in
      N_{B_{2}}(c_{2})$, then 
      $F = c_{1}Tc_{1}a_{1}P_{H}^{1}a_{2}c_{2}'T''c_{2}P_{H}^{2}c_{1}$ is
      a $[2, 2s]$-factor in $G^{2}$, contradicting condition (2).

 $b)$ The graph $B_{1}^{2}$ has a spanning $s$-trail $T^{*}$ between $c_{1}$ and some neighbor $c_{1}' \in
      N_{B_{1}}(c_{1})$.

      If the graph $B_{2}^{2}$ has a spanning closed $s$-trail $T^{**}$ in which $d_{T^{**}}(c_{2}) \leq 2s - 2$, then
      $F = c_{1}T^{*}c_{1}'a_{1}P_{H}^{1}a_{2}c_{2}T^{**}c_{2}P_{H}^{2}c_{1}$ is a $[2, 2s]$-factor in $G^{2}$, a contradiction.

      If the graph $B_{2}^{2}$ has a spanning $s$-trail $T^{\bullet}$ between $c_{2}$ and some neighbor
      $c_{2}' \in N_{B_{2}}(c_{2})$, then 
      $F = c_{1}T^{*}c_{1}'a_{1}P_{H}^{1}a_{2}
      c_{2}'T^{\bullet}c_{2}P_{H}^{2}c_{1}$ is a $[2, 2s]$-factor in $G^{2}$ and contradicts (2).  \hfill $\blacksquare$

The graph $G$ in Figure 1 shows that (for $s=1$) the constant $3$ in Theorem \ref{thinduced} cannot be decreased. Although every induced $S(K_{1,2s+1})$ in $G$ has at least two edges in a  block of degree at most two, $G^2$ has no $[2,2s]$-factor.

\begin{figure}
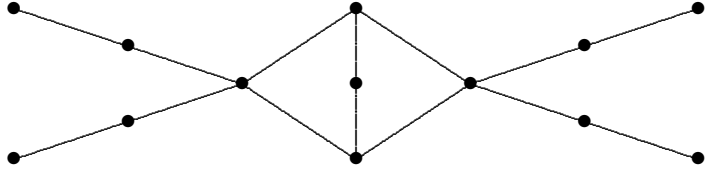

 \beginpicture
\setcoordinatesystem units <1mm,1mm>
\setplotarea x from 0 to 155, y from 15 to 40
\large 
\put{$\bullet$} at 30 40
\put{$\bullet$} at 30 20
\put{$\bullet$} at 45 35
\put{$\bullet$} at 45 25
\put{$\bullet$} at 60 30
\put{$\bullet$} at 75 40
\put{$\bullet$} at 75 20
\put{$\bullet$} at 90 30
\put{$\bullet$} at 105 35
\put{$\bullet$} at 105 25
\put{$\bullet$} at 120 40
\put{$\bullet$} at 120 20
\put{$\bullet$} at 75 30
\plot 30 40  60 30  30 20 /
\plot 120 40  90 30  120 20 /
\plot 75 40  75 20  90 30  75 40  60 30  75 20 /
\endpicture
\caption{An example showing that a condition in Theorem 2 cannot be relaxed.}
\end{figure}

\end{document}